\documentclass[a4paper]{jpconf}
\usepackage{graphicx}
\begin{document}

\title{REAL-TIME EQUILIBRIUM RECONSTRUCTION IN A TOKAMAK}

\author{J. Blum, C. Boulbe and B. Faugeras}

\address{Laboratoire J.A. Dieudonn\'e, UMR 6621, Universit\'e de Nice
  Sophia Antipolis, Parc Valrose, 06108 Nice Cedex 02, France}

\ead{jblum@unice.fr, boulbe@unice.fr, faugeras@unice.fr}

\begin{abstract}
This paper deals with the numerical reconstruction 
of the plasma current density in a Tokamak and of its equilibrium. 
The problem consists in the identification of a non-linear source 
in the 2D Grad-Shafranov equation, which governs the axisymmetric equilibrium 
of a plasma in a Tokamak. 
The experimental measurements that enable this identification 
are the magnetics on the vacuum vessel, 
but also polarimetric and interferometric measures on several chords, 
as well as motional Stark effect or pressure measurements. 
The reconstruction can be obtained in real-time using 
a finite element method, a non-linear fixed-point algorithm 
and a least-square optimization procedure.
\end{abstract}



\section{Introduction}

The problem of the equilibrium of a plasma in a Tokamak is 
a free boundary problem in which the plasma boundary is defined either 
by its contact with a limiter or as being a magnetic separatrix. 
Inside the plasma, the equilibrium equation
in an axisymmetric configuration is called Grad-Shafranov equation 
\cite{Grad:1958,Shafranov:1958,Mercier:1974}. 
The right-hand side of this equation is a non-linear source 
which represents the toroidal component of the plasma current density.

An important problem is the identification of this non-linearity 
\cite{Lao:1990,Blum:1990,Blum:1997}. The aim of this paper is to present a method for 
real-time identification from experimental measurements, such as magnetic measurements on the vacuum vessel, polarimetric measurements (integrals of the magnetic 
field over several chords), MSE (Motional Stark Effect) 
and pressure measurements. 
The pressure is supposed to be isotropic. For the anisotropic pressure case, 
one can refer to \cite{Zwingmann:2000}.

The next section is devoted to the mathematical modelling 
of the equilibrium problem in axisymmetric configurations. 
The inverse reconstruction problem is adressed in the last section.

\section{Mathematical modelling of axisymmetric equilibrium 
of the plasma in a Tokamak}

The equations which govern the equilibrium of a plasma in the presence of a magnetic field are on the one hand Maxwell's equations and on the other hand the equilibrium equations for the plasma itself.

The magnetostatic Maxwell's equations as follows are satisfied in the whole of space (including the plasma):

\begin{equation}
\label{eqn:maxwell}
\left \lbrace
\begin{array}{lll}
\nabla \cdot B &=& 0\\
\nabla \times (\displaystyle \frac{B}{\mu}) &=&  j 
\end{array}
\right.
\end{equation}

where $B$ represents the magnetic field, $\mu$ is the magnetic permeability and 
$j$ is the current density. 
The first relation of (\ref{eqn:maxwell}) is the equation of conservation of magnetic induction and the second one is Ampere's Theorem.

The momentum equation for a plasma is
\begin{equation}
\label{eqn:momentum}
\rho \frac{du}{dt}=j \times B - \nabla p
\end{equation}
where $u$ represents the mean velocity of particles 
and $\rho$ the mass density.
At the resistive time-scale the first term can be neglected \cite{Grad:1970} and the equilibrium equation for the plasma is 
\begin{equation}
\label{eqn:equilibrium}
\nabla p = j \times B
\end{equation}
This equation (\ref{eqn:equilibrium}) means that the plasma is in equilibrium when the force $\nabla p$ due the kinetic pressure $p$ is equal to the Lorentz force of the magnetic pressure $j \times B$. We deduce immediately from (\ref{eqn:equilibrium}) that
\begin{equation}
\label{eqn:cons1}
B \cdot \nabla p = 0
\end{equation}
\begin{equation}
\label{eqn:cons2}
j \cdot \nabla p = 0
\end{equation}
Thus for a plasma in equilibrium the field lines and the current lines lie on 
isobaric surfaces ($p=const.$); these surfaces, generated by the field lines, 
are called magnetic surfaces. In order for them to remain within a bounded volume 
of space it is necessary that they have a toroidal topology. These surfaces form a family 
of nested tori. 
The innermost torus degenerates into a curve which is called the magnetic axis.

In a cylindrical coordinate system $(r,z,\phi)$ (where $r=0$ is the major axis of the torus) 
the hypothesis of axial symmetry consists in assuming that the 
magnetic field $B$ is independent of the toroidal angle $\phi$. 
The magnetic field can be decomposed as $B=B_p +B_{\phi}$, where $B_p=(B_r,B_z)$ is 
the poloidal component and $B_{\phi}$ is the toroidal component. 
From equation (\ref{eqn:maxwell}) one can define the poloidal flux $\psi(r,z)$ such that 
\begin{equation}
\left \lbrace
\begin{array}{lll}
B_r &=&-\displaystyle \frac{1}{r}\frac{\partial \psi}{\partial z}\\[10pt]
B_z &= & \displaystyle \frac{1}{r}\frac{\partial \psi}{\partial r} 
\end{array}
\right.
\end{equation}
Concerning the toroidal component $B_{\phi}$ we define $f$ by 
\begin{equation}
B_{\phi}=\frac{f}{r}e_{\phi}
\end{equation}
where $e_{\phi}$ is the unit vector in the toroidal direction, and $f$ is the 
diamagnetic function. 
The magnetic field can be written as:
\begin{equation}
\label{eqn:B}
\left \lbrace
\begin{array}{lll}
B&=&B_p+B_{\phi}\\[10pt]
B_p &=&\displaystyle \frac{1}{r}[\nabla \psi \times e_{\phi}] \\[10pt]
B_{\phi} &= & \displaystyle \frac{f}{r} e_{\phi} 
\end{array}
\right.
\end{equation}
According to (\ref{eqn:B}), in an axisymmetric configuration the magnetic surfaces 
are generated by the rotation of the flux lines $\psi=const.$ around the axis $r=0$ of the torus.

From (\ref{eqn:B}) and the second relation of (\ref{eqn:maxwell}) we obtain the 
following expression for $j$:
\begin{equation}
\label{eqn:j}
\left \lbrace
\begin{array}{lll}
j&=&j_p+j_{\phi}\\[10pt]
j_p &=&\displaystyle \frac{1}{r}[\nabla (\frac{f}{\mu}) \times e_{\phi}] \\[10pt]
j_{\phi} &= & (-\Delta^* \psi) e_{\phi} 
\end{array}
\right.
\end{equation}
where $j_p$ and $j_{\phi}$ are the poloidal and toroidal components respectively of $j$, and the operator $\Delta^*$ is defined by
\begin{equation}
\label{eqn:delta*}
\Delta^* .=  \frac{\partial }{\partial r}(\frac{1}{\mu r}  \frac{ \partial . }{\partial r}) 
+ \frac{\partial}{\partial z}(\frac{1}{\mu r} \frac{ \partial .}{\partial z})
\end{equation}

Expressions (\ref{eqn:B}) and (\ref{eqn:j}) for $B$ and $j$ are valid in the whole of space since they involve only Maxwell's equations and the hypothesis of axisymmetry.

In the plasma region, relation (\ref{eqn:cons1}) implies that $\nabla p$ and 
$\nabla \psi$ are colinear, and therefore $p$ is constant on each magnetic surface. This can be denoted by
\begin{equation}
p=p(\psi)
\end{equation}
Relation (\ref{eqn:cons2}) combined with the expression (\ref{eqn:j}) implies 
that $\nabla f$ and $\nabla p$ are colinear, and therefore $f$ is likewise constant on each magnetic surface
\begin{equation}
f=f(\psi)
\end{equation}
The equilibrium relation (\ref{eqn:equilibrium}) combined with the expression 
(\ref{eqn:B}) and (\ref{eqn:j}) for $B$ and $j$ implies that:
\begin{equation}
\label{eqn:gradshaf0}
\nabla p = - \displaystyle \frac{\Delta^* \psi}{r} \nabla \psi - \frac{f}{\mu_0 r^2} \nabla f
\end{equation}
which leads to the so-called Grad-Shafranov equilibrium equation:
\begin{equation}
\label{eqn:gradshaf}
-\Delta^* \psi = r p'(\psi) + \frac{1}{\mu_0 r}(ff')(\psi)
\end{equation}
where $\Delta^*$ is the linear elliptic operator given by (\ref{eqn:delta*}) 
in which $\mu$ is equal to the magnetic permeability $\mu_0$ of the vacuum. 

From (\ref{eqn:j}) it is clear that right-hand side of (\ref{eqn:gradshaf}) 
represents the toroidal component of the plasma current density. It involves 
functions $p(\psi)$ and $f(\psi)$ which are not directly measured inside the plasma.

\begin{figure}[!h]
\centering
\includegraphics[width=9cm,height=11cm,angle=-90]{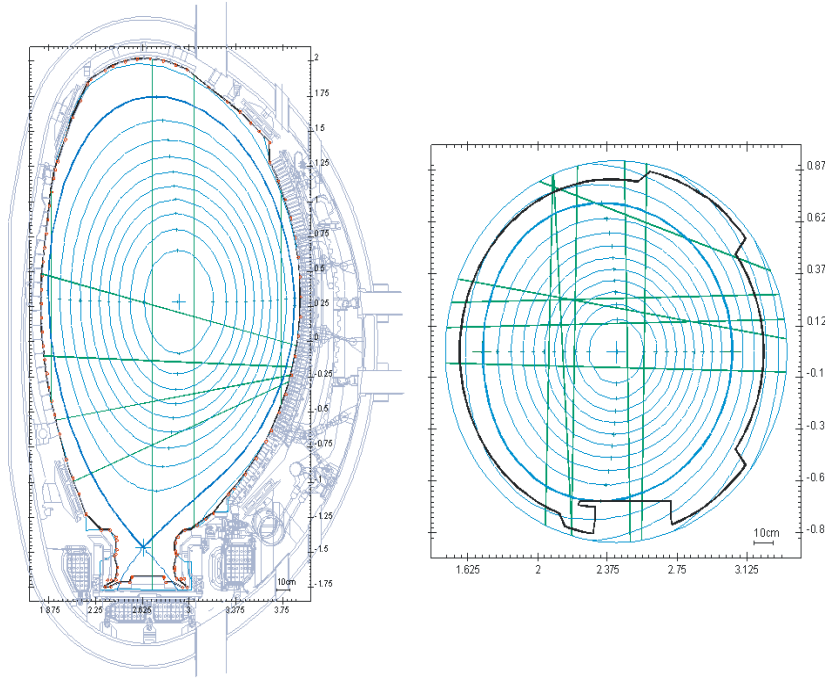} 
\caption{\label{fig:jet-tore} Definition of the plasma boundary (thick blue line). 
Left, JET (Joint European Torus) example, X-point configuration. Right, TORE SUPRA (the CEA-EURATOM Tokamak at Cadarache) example, limiter configuration (the limiter is represented by the black line). The thin blue lines represent iso-contours of $\psi$.}
\end{figure}

In the vacuum, the magnetic flux $\psi$ satisfies
\begin{equation}
\label{eqn:gradshafvac}
-\Delta^* \psi = 0
\end{equation}
The equilibrium of a plasma in a domain $\Omega$ representing the vacuum 
region is a free boundary problem. 
The plasma free boundary is defined either by its contact with a limiter $D$ 
(outermost flux line inside the limiter) 
or as being a magnetic separatrix (hyperbolic line with an $X$-point, $X$). 
The region $\Omega_p \subset \Omega$ containing the plasma is defined as
\begin{equation}
\Omega_p=\lbrace x \in \Omega ,\ \psi(x) \ge \psi_b \rbrace
\end{equation}
where either $\psi_b=\displaystyle \max_{D} \psi$ in the limiter configuration 
or $\psi_b=\psi(X)$ in the $X$-point configuration (see Fig. \ref{fig:jet-tore})

Assuming Dirichlet boundary conditions, $h$, are given on 
$\Gamma = \partial \Omega$ which is the poloidal 
cross-section of the vacuum vessel, 
the final equations governing the behaviour 
of $\psi(r,z)$ inside the vacuum vessel, are:
\begin{equation}
\label{eqn:final}
\left \lbrace
\begin{array}{rcl}
-\Delta^* \psi &= &[r A(\bar{\psi}) +  \displaystyle 
\frac{1}{r} B(\bar{\psi})] \chi_{\Omega_p}\quad \mathrm{in}\ \Omega \\[10pt]
\psi&=& h\quad  \mathrm{on}\  \Gamma 
\end{array} 
\right.
\end{equation}
with 
$A(\bar{\psi})= p'(\bar{\psi})$ and 
$B(\bar{\psi})=\displaystyle \frac{1}{\mu_0}(ff')(\bar{\psi})$, 
$\bar{\psi}= \displaystyle \frac{\psi - \displaystyle \max_\Omega \psi}{\psi_b - \displaystyle \max_\Omega \psi} \in [0,1]$ 
in $\Omega_p$ 
(this normalized flux is introduced so that $A$ and $B$ 
are defined on the fixed interval $[0,1]$), $\chi_{\Omega_p}$ 
is the characteristic function of $\Omega_p$.

The aim of the following section of this paper is to provide a method 
for the real-time identification of the plasma current i.e. 
the non-linear functions $A$ and $B$ in the 
elliptic equation (\ref{eqn:final}).

\section{The inverse problem}

\subsection{Experimental measurements}
\begin{figure}[!h]
\centering
\includegraphics[width=10cm,height=12cm,angle=-90]{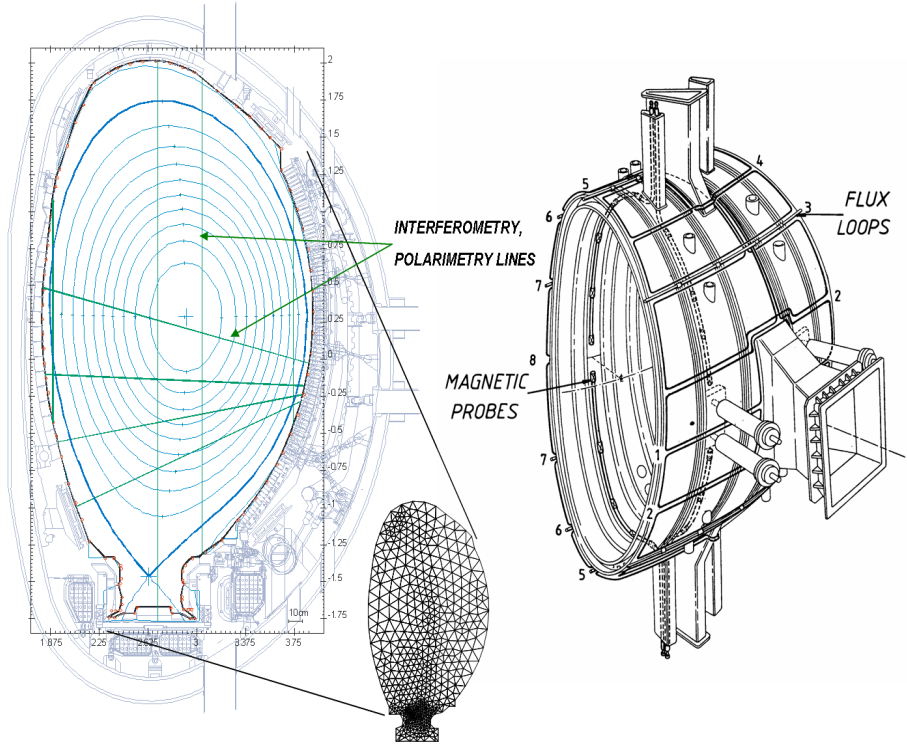} 
\caption{\label{fig:jet-data-mesh} Left: the straight green lines represents the chords used 
for polarimetry and interferometry measurements. Right: part of the vacuum vessel. At the bottom middle an example of finite element mesh used for numerical simulations (see next Section).}
\end{figure}

The given experimental measurements are:  

\begin{itemize}

\item the magnetic measurements
 \begin{itemize}
   \item $\psi(M_i)=h_i\  \mathrm{on}\  \Gamma$, 
   given by the flux loops (see Fig. \ref{fig:jet-data-mesh}). Thanks to an interpolation between points $M_i$ these measurements provide the Dirichlet boundary condition $h$. 
   \item $\displaystyle \frac{1}{r} \frac{\partial \psi}{\partial n}(N_i)=g_i\ \mathrm{on}\ \Gamma$, 
   which corresponds to the component of the magnetic poloidal field, 
   measured by the magnetic probes (see Fig. \ref{fig:jet-data-mesh}), 
   which is tangent to the vacuum vessel. Indeed from Eq. (\ref{eqn:B}) the tangential 
   component of $B_p$ 
   is equal to the normal component $\displaystyle \frac{1}{r}\frac{\partial \psi}{\partial n}$ of $\displaystyle \frac{1}{r}\nabla \psi$.
 \end{itemize}

\item the polarimetric measurements which give the Faraday rotation of the angle of infrared radiation crossing the section of the plasma along several chords $C_i$:
$$
\displaystyle \int_{C_i} n_e(\bar{\psi}) B_{\|}dl = 
\displaystyle \int_{C_i} \frac{n_e(\bar{\psi})}{r} \frac{\partial \psi}{\partial n}dl=\alpha_i
$$
where $n_e$ represents the electronic density which is approximately constant on each flux line, 
$B_{\|}$ is the component of the poloidal field tangent to $C_i$ 
and $\displaystyle \frac{\partial}{\partial n}$ represents the normal 
derivative of $\psi$ with respect to $C_i$. 

\item the interferometric measurements which give the density integrals over the chords $C_i$
$$
\displaystyle \int_{C_i}n_e(\bar{\psi})dl=\beta_i 
$$

\item the kinetic pressure measurements obtained from density and temperature measurements, for instance in the equatorial plane:
$$
p(r,0)=p_d(r)
$$

\item and MSE (Motional Stark Effect) angle measurements 
taken at different points $x_i=(r_i,z_i)$:
$$
m(B_r,B_z,B_\phi)_i=\gamma_i
$$
with
$$
\tan ( m(B_r,B_z,B_\phi) ) = \displaystyle \frac{a_{1} B_r + a_{2} B_z + a_{3} B_\phi}
{a_{4} B_r + a_{5} B_z + a_{6} B_\phi}
$$
\end{itemize}

\subsection{Statement of the inverse problem}

The numerical identification problem is formulated as a least-square 
minimization with a Tikhonov regularization. 
The cost function is defined as:

\begin{equation}
\label{eqn:costfunction}
J(A,B,n_e)=J_0 + K_1 J_1 + K_2 J_2 + K_3 J_3 + K_4 J_4 + J_\epsilon
\end{equation}

with

$$
\begin{array}{l}
 J_0=  \displaystyle \sum_i 
 ( \displaystyle \frac{1}{r} \frac{\partial \psi}{\partial n}(N_i) 
 - g_i)^2  \\[10pt]
 J_1=  \displaystyle \sum_i (\displaystyle \int_{C_i} 
 \frac{n_e}{r} \frac{\partial \psi}{\partial n}dl - \alpha_i)^2  \\[10pt]
J_2= \displaystyle \sum_i (\displaystyle \int_{C_i} n_edl - \beta_i)^2 \\[10pt]
J_3=  \displaystyle \int_{R_{min}}^{R_{max}} (p(r,0)-p_d(r))^2 dr\\[10pt]
J_4=  \displaystyle \sum_i (  
m(B_r,B_z,B_\phi)_i - \gamma_i )^2
\end{array}
$$

and $K_1$, $K_2$, $K_3$ and $K_4$ are weighting parameters enabling to give more or less importance to the corresponding experimental measurements \cite{Blum:1990}.  

The inverse problem of the determination of $A$ and $B$ is ill-posed. 
Hence a regularization procedure can be used to transform it into a well-posed one 
\cite{Tikhonov:1977}. 
The Tikhonov regularization term $J_\epsilon$ constrains 
the function $A$, $B$ and $n_e$ to be smooth enough and reads:
$$
J_\epsilon = 
\epsilon_1 \displaystyle \int_0^1 [A''(x)]^2 dx 
+
\epsilon_2 \displaystyle \int_0^1 [B''(x)]^2 dx 
+
\epsilon_3 \displaystyle \int_0^1 [n_e''(x)]^2 dx 
$$
where $\epsilon_1$, $\epsilon_2$ and $\epsilon_3$ are the regularizing parameters. 

It should be noticed that the electronic density $n_e$ does not intervene in 
Eq. (\ref{eqn:final}). However as soon as we want to use the polarimetric measurements it is necessary to include $n_e$ (and hence interferometry) in the identification procedure.
The inverse problem can finally be formulated as,

\begin{equation}
\left \lbrace
\begin{array}{l}
\mathrm{Find}\ A^*,\ B^*,\ n_e^*\ \mathrm{such}\ \mathrm{that}: \\[10pt]
J(A^*,B^*,n_{e}^*)=\inf J(A,B,n_{e}) 
\end{array}
\right.
\end{equation}

\subsection{Numerical identification}

Problem (\ref{eqn:final}) is solved using a finite element method 
\cite{Ciarlet:1980}.
Let $H^1(\Omega)$ and $V=H^1_0(\Omega)$ denote the usual Sobolev spaces. 
The finite element approximation is based on the following weak formulation:
\begin{equation}
\label{eqn:var}
\left \lbrace
\begin{array}{l}
\mathrm{Find}\ \psi \in H^1(\Omega),\ \mathrm{such}\ \mathrm{that}\ 
\psi=h\ \mathrm{on}\ \Gamma,\ \mathrm{and} \\[10pt]
\displaystyle \int_\Omega \displaystyle \frac{1}{\mu_0 r}\nabla \psi \cdot \nabla v dx = 
\int_{\Omega_p}[r A(\bar{\psi})+ \frac{1}{r}B(\bar{\psi})] v dx\quad \forall v\ \in V \\[10pt]
\end{array}
\right .
\end{equation}
Classically $\Omega$ is approximated using triangles by a polygonal domain $\Omega_h$, 
the space $V$ is approximated by a space $V_h$ of finite dimension $n$. A $P1$ finite element 
method is used, in which functions of $V_h$ are affine over each triangle and continuous on the 
whole domain.

Let $K$ denote the finite element stiffness matrix. 
Let us also (abusively) denote by $\psi \in \mathbf{R}^n$ the components of 
the magnetic flux fonction approximated in $V_h$.

The unknown functions $A$, $B$ and $n_e$ are approximated by a decomposition 
in a reduced basis $(\phi_i)_{i=1,\dots m}$ 
$$
\begin{array}{l}
A(x)=\displaystyle \sum_i a_i \phi_i(x)\\[10pt]
B(x)=\displaystyle \sum_i b_i \phi_i(x)\\[10pt]
n_e(x)=\displaystyle \sum_i c_i \phi_i(x)
\end{array}
$$
This basis can be made of different types of functions 
(polynomials, splines, wavelets, etc \dots) \cite{Blum:1997}. 
Let $u$ be the vector of $\mathbf{R}^{3m}$ defined by 
$u=(a_1,\dots ,a_m,b_1,\dots,b_m,c_1,\dots,c_m)$.
With these notations the discretization of problem (\ref{eqn:var}) 
can be written as follows:

\begin{equation}
\label{eqn:fp}
\left \lbrace
\begin{array}{l}
\mathrm{Given}\ u \in \mathbf{R}^{3m},\ \mathrm{solve}\ \mathrm{the}\ \mathrm{fixed-point}\ \mathrm{equation}\\[10pt]
\tilde{K} \psi = D(\psi) u + h
\end{array}
\right .
\end{equation}

Where $D(\psi)$ denotes the $n\times 3m$ ``plasma current matrix'', 
and $\tilde{K}$ is the stiffness matrix modified in order to impose 
the Dirichlet boundary condition represented by $h$.


The discrete inverse optimization problem is:

\begin{equation}
\label{eqn:optim}
\left \lbrace
\begin{array}{l}
\mathrm{Find}\ u\ \mathrm{minimizing}:\\[10pt]
J(u)=\|  C(\psi)\psi - k \|^2 + u^{T} \Lambda u\\[10pt] 
\mathrm{with}\ \psi\ \mathrm{satisfying}\ \mathrm{(\ref{eqn:fp})} 
\end{array}
\right.
\end{equation}

where $C(\psi)$ is the observation operator. 
The quantity $C(\psi)\psi$ represents the outputs of the model 
corresponding to the experimental measurements, given in a previous subsection, denoted by $k$. 
The matrix $\Lambda$ represents the regularization terms. 
The first term of $J$ in Eq. (\ref{eqn:optim}) corresponds to 
$J_0 + K_1 J_1 + K_2 J_2 + K_3 J_3 + K_4 J_4$ and the second to $J_\epsilon$.

In order to solve this problem we use 
an iterative algorithm based on fixed-point iterations 
for Eq. (\ref{eqn:fp}) and the normal equation of Eq. (\ref{eqn:optim}).

\subsubsection{Algorithm}
At the $n$-th iteration, $\psi_n$ and $u_n$ are given. 
The non-linear mapping $u \mapsto \psi(u)$ is approximated by the affine relation
$$
\psi =\tilde{K}^{-1}[D(\psi_n) u +h]
$$
and the cost function to be minimized by
$$
\begin{array}{lll}
J(u)&=&\|  C(\psi_n)\psi - k \|^2 + u^{T} \Lambda u \\[10pt]
&=&\|  C(\psi_n) \tilde{K}^{-1} D(\psi_n) u +( C(\psi_n) \tilde{K}^{-1} h - k ) \|^2 + u^{T} \Lambda u \\[10pt]
&=&\|  E_n u+ F_n \|^2 + u^{T} \Lambda u 
\end{array}
$$
with obvious notations. 
The normal equation 
$$
(E^{T}_n E_n + \Lambda)u= - E^{T}_n F_n
$$
is solved to update $u_n$ to $u_{n+1}$. 
Then a fixed-point iteration for Eq. (\ref{eqn:fp}) 
enables the update of $\psi_n$
to $\psi_{n+1}$
$$
\psi_{n+1}= \tilde{K}^{-1}[D(\psi_n) u_{n+1} +h].
$$
Since the algorithm is usually initialized with the equilibrium at 
a previous time step, two or three fixed-point iterations are usually 
enough to ensure convergence.


\subsubsection{Equinox software}
\begin{figure}[!h]
\centering
\includegraphics[width=10cm,height=12cm,angle=-90]{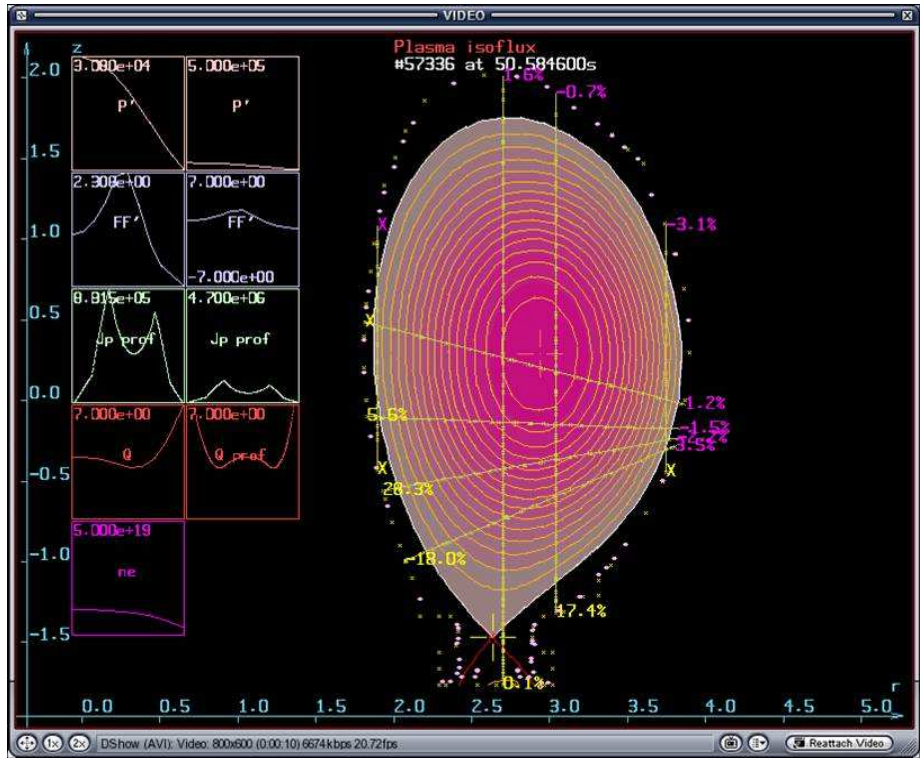} 
\caption{\label{fig:equinox} An output of EQUINOX. The plasma is in an X-point configuration. 
On the left column the identified $p'$, $ff'$ and $n_e$ functions as well as 
the toroidal current density $j$ and the safety factor $q$ are displayed 
in terms of $\psi$ and of $r$ (in the equatorial plane).}
\end{figure}

Based on the algorithm presented above, a C++ software, 
called EQUINOX \cite{Bosak:2001,Blum:2004,Bosak:2003} has been developed 
in collaboration with the Fusion Department at Cadarache, 
and has been implemented for JET (Joint European Torus) 
and for TORE SUPRA (the CEA-EURATOM Tokamak at Cadarache). 
Figure \ref{fig:equinox} shows a graphical output of Equinox. 
With all these techniques it is possible to follow the quasi-static 
evolution of the plasma equilibrium, either in TORE SUPRA or JET configurations, 
with free boundaries defined either by limiter contact or with an X-point. 
It is also possible to simulate ITER configurations.

\section{Conclusion}

We have presented an algorithm for the identification of the current 
density profile in Grad-Shafranov equation from experimental measurements.\\
The decomposition of the unknown functions $p'(\psi)$ and $ff'(\psi)$ 
in a reduced basis makes it possible to do the reconstruction in real-time.\\
The choice of this reduced basis must still be improved and optimized 
(robustness, precision, ...).\\
Real-time reconstruction makes possible 
future real-time control of the current profile \cite{Joffrin:2003}.

\section{Aknowledgements}
The authors are very thankful to K. Bosak for his enormous work in 
developing the real-time EQUINOX software and to E. Joffrin and S. Bremond for 
their 
contributions into the domain of real-time identification and control.

\section*{References}
\bibliographystyle{iopart-num}

\bibliography{biblio-fusion}

\end{document}